# BEST SUBSET SELECTION, PERSISTENCE IN HIGH-DIMENSIONAL STATISTICAL LEARNING AND OPTIMIZATION UNDER $L_1$ CONSTRAINT

BY EITAN GREENSHTEIN

*Purdue University*

Let $(Y, X_1, \ldots, X_m)$ be a random vector. It is desired to predict $Y$ based on $(X_1, \ldots, X_m)$. Examples of prediction methods are regression, classification using logistic regression or separating hyperplanes, and so on.

We consider the problem of best subset selection, and study it in the context $m = n^\alpha$, $\alpha > 1$, where $n$ is the number of observations. We investigate procedures that are based on empirical risk minimization. It is shown, that in common cases, we should aim to find the best subset among those of size which is of order $o(n/\log(n))$. It is also shown, that in some "asymptotic sense," when assuming a certain sparsity condition, there is no loss in letting $m$ be much larger than $n$, for example, $m = n^\alpha$, $\alpha > 1$. This is in comparison to starting with the "best" subset of size smaller than $n$ and regardless of the value of $\alpha$.

We then study conditions under which empirical risk minimization subject to $l_1$ constraint yields nearly the best subset. These results extend some recent results obtained by Greenshtein and Ritov.

Finally we present a high-dimensional simulation study of a "boosting type" classification procedure.

**1. Introduction and preliminaries.** Let $Z^i = (Y^i, X_1^i, \ldots, X_m^i)$, $i = 1, \ldots, n$, be i.i.d. vectors, $Z^i \sim F$ where $F$ is unknown. It is desired to find a good predictor for $Y$ given $X_1, \ldots, X_m$, based on the observations $Z^i, i = 1, \ldots, n$. In this paper we consider high-dimensional learning problems, where the objective is to select a good predictor from a large class, based on minimizing an empirical risk. We concentrate on the case where the dimension is much larger than the number of observations, that is, $m \gg n$.

There are three main goals of this paper. One is to advocate the practice of turning to high dimensions of explanatory variables for the purpose of finding good predictors. Another is to give a perspective to the phenomenon of "not getting overfit," when applying high-dimensional procedures, as discussed in [2]. We will suggest that often such procedures may be viewed as (suboptimal) optimization methods for finding the empirically best subset of explanatory variables. A final goal is to show that often optimization under $l_1$ constraint (as in "Lasso") could be a helpful and computationally feasible method for finding good predictors in high dimensions.

We describe now a few examples where an analysis with $m \gg n$ is conducted. In microarray experiments the explanatory variables are measurements describing activity of certain $m$ genes in $n$ subjects, while the response could be survival time or an indicator of the event that the subject has a certain disease, and so on; see [21]. Under the current technology, a typical microarray experiment involves thousands of genes, that is, the dimension $m$ is of the order of thousands, while $n$ is of the order of hundreds or less.

In [25], page 496, the following pattern recognition example is described. It is desired to train a machine to identify handwritten digits for the purpose of recognizing handwritten zip codes. The raw data given to the machine comes from 256 pixels, that is, the raw data is made up of 256 variables. Yet, for their classification method, they considered all interactions up to order 7. This creates $m \approx 10^{16}$ explanatory variables constructed from the initial set of 256. The amount of data (or training set) they were using was $n = 7291$.

Finally, consider the following example as a plausible data mining application of analysis with $m \gg n$. An insurance company is interested in estimating the probability of a claim, due to a car accident, by various customers. We may define for each customer quite a few categorical variables based on age, sex, car make, car model, marital status, address, and so on. Considering also third- or fourth-order interactions of these categorical variables, one does not need a lot of imagination to come up with tens and hundreds of millions of categorical explanatory variables. Of course, the insurance company might have access to a big historical database, so $n$ may also be very large.

Although our motivation is to understand the problem where $m \gg n$, there are also implications to the following more classical problem when $m < n$. Informally the problem may be stated as follows: how many observations, $n$, do we need, in order to accurately estimate $m$ parameters? Our asymptotic approach suggests that in many cases the condition $m \log(m) = o(n)$ suffices. See further discussion at the end of this section.

We will consider and formulate our problem in various degrees of generality. The ideas are easier to introduce and motivate through the problem of best subset selection in regression, but will be carried out in a more general context.

Let $Z = (Y, X_1, \ldots, X_m)$ be a random vector $Z \sim F$, $F$ unknown. Consider first the problem of selecting a linear predictor for $Y$ based on $X_1, \ldots, X_m$, that is, a function of the form $\sum_j \beta_j X_j$. We identify a predictor with the vector $\beta = (\beta_1, \ldots, \beta_m)$. Its performance is evaluated based on

$$(1) \qquad L_F(\beta) = E_F\left(Y - \sum \beta_j X_j\right)^2.$$

The selection of a predictor is based on a sample of i.i.d. observations $Z^i$, $i = 1, \ldots, n$. In practice, as the sample size, $n$, increases, we might want to consider more complicated models or linear predictors, that is, increase the number $m$ of explanatory variables. Thus, a worthwhile asymptotic study is of a triangular array form, where we are given $n$ i.i.d. observations $Z_n^1, \ldots, Z_n^n$ at stage $n$, $Z_n^i \sim F_n$, $F_n$ is unknown, $F_n \in \mathcal{F}_n$. In order to simplify notation, we will drop the index $n$ of the triangular array and write $Z^i$; thus, at stage $n$, $Z^i = (Y^i, X_1^i, \ldots, X_m^i)$. Here $m = m(n)$ is the number of explanatory variables, which depends on $n$ and typically grows with $n$. We will study asymptotics where $m = n^\alpha$, $\alpha > 1$. See further discussion on the triangular array setup in [13]. Further papers investigating a similar regression triangular array structure are [16, 17, 19]. These papers also study the Lasso and regularization via $l_1$ constraints, as we do in this paper. A recent paper that studies the virtue of letting $m$ be much larger than $n$ in classification problems is [1].

The above regression setup motivates us to generalize as follows. Consider a triangular array as before, equipped with an abstract triangular structure of parametrized predictors, that is, at stage $n$, a collection of functions

$$\{g_\beta, \beta \in B^n\},$$

where $g_\beta = g_\beta(X_1, \ldots, X_{m(n)})$, and the parametrization is Euclidean. Consider a general nonnegative prediction loss $l$, incurred for predicting $g_\beta(X_1, \ldots, X_{m(n)})$ when the outcome is $Y$,

$$l = l(Y, g_\beta(X_1, \ldots, X_{m(n)})).$$

To simplify notation, we will abuse and write

$$l(\beta, Z) \equiv l(Y, g_\beta(X_1, \ldots, X_{m(n)})).$$

As in equation (1), we define

$$(2) \qquad L_F(\beta) = E_F l(\beta, Z).$$

Note in equation (1) we used a squared loss $l$. As an additional example, consider classification where $Y$ may be either $+1$ or $-1$, the predictors are of the type $g_\beta(X_1, \ldots, X_m) = \text{sign}(\sum \beta_j X_j)$, and the prediction loss is 0–1.

In the current more abstract formulation, we will consider entry $j$ of the parameter $\beta$ "*active*" if $\beta_j \neq 0$. Note, in order to relate to regression and other important examples, we denote both the dimension of the explanatory variables and of the parameter space by $m$. However, in the abstract formulation the dimension of the explanatory variables is actually not relevant. In the sequel, assumptions made about $m = m(n)$ may in fact be assumed only on the dimension of the parameter space.

Let
$$\beta^*_{F_n} = \arg\min_{\beta \in B^n} L_{F_n}(\beta).$$

From now, when we say *triangular array*, we mean a sequence of collections of distributions $\mathcal{F}_n$, a sequence $B^n$ of collections of predictors which are available at stage $n$, $n = 1, 2, \ldots$, and a prediction loss function $l$.

We will study sequences of procedures $\hat{\beta} = \hat{\beta}(Z^1, \ldots, Z^n)$ that select a predictor $\beta \in B^n$, based on the observations $Z^1, \ldots, Z^n$. Here $Z^i$ are i.i.d. distributed $F_n$, $F_n \in \mathcal{F}_n$. The dependence of $\hat{\beta} \equiv \hat{\beta}_n$ on $n$ is often suppressed, and we will loosely say the procedure $\hat{\beta}$.

DEFINITION 1. Given a triangular array, the sequence of procedures $\hat{\beta}_n$ is persistent with respect to $B^n$ if, for every $\varepsilon > 0$,

(3) $$\sup_{F_n \in \mathcal{F}_n} P_{F_n}(L_{F_n}(\hat{\beta}_n) - L_{F_n}(\beta^*_{F_n}) > \varepsilon) \to 0.$$

It is not difficult to see that the above is equivalent to the following: for any sequence $F_n \in \mathcal{F}_n$,
$$L_{F_n}(\hat{\beta}_n) - L_{F_n}(\beta^*_{F_n}) \xrightarrow{p} 0.$$

Here the distribution of $\hat{\beta}_n$ is determined by $F_n$.

REMARK 1. (a) The concept of persistence is close to that of consistency. Yet, in consistency there is a certain, usually "true," fixed parameter to which a consistent estimator converges. In our setup the analog of the true parameter is $\beta^*_{F_n}$, which changes with $n$. Also, in consistency convergence is usually in terms of the Euclidean distance between the true parameter and its estimator, while in persistence the distance is tied to the loss.

(b) Consider the triangular array structure that motivates us, where as $n$ grows we consider larger nested collections of predictors $B^n$. In such a nested structure we may consider the joint distribution $F^0_\infty$ of all variables, that is, the joint distribution of $(Y, X_1, \ldots, X_{m(\infty)})$. Let $F^0_n$ be the marginal of $F^0_\infty$ on $\sigma(Y, X_1, \ldots, X_{m(n)})$. Obviously $L_{F^0_n}(\beta^*_{F^0_n})$ is monotone decreasing since $B^n \subset B^{n+1}$. Thus, there is a limit
$$\lim_n L_{F^0_n}(\beta^*_{F^0_n}) = r(F^0_\infty).$$

When $r(F^0_\infty) > 0$, the persistence criterion should have appeal. In situations where $r(F^0_\infty) = 0$, other criteria should be studied and rates of convergence become relevant, rather than only persistence.

Under mild conditions, existence of a persistent procedure will follow if $\beta^*_{F_n}$, the best predictor in $B^n$, has $k_n = o(n/\log(n))$ nonzero entries, also termed $o(n/\log(n))$ *sparsity rate*. This may be shown by a simple entropy derivation; see Theorem 1 of the next section. It will also be demonstrated in Section 2 that if, for the relevant sequence $B^n$, the corresponding sequence $\beta^*_{F_n}$ has $o(n/\log(n))$ sparsity rate, then there is only a mild effect on the ability to find a predictor which is nearly as good as $\beta^*_{F_n}$, when increasing the dimension $m$ dramatically.

*Discussion of the asymptotics and the sets $B^n$.* We further discuss now our notion of persistence with respect to sets $B^n$. The discussion is in light of the regression setup with $m \gg n$. Usually in asymptotics we evaluate procedures comparing their estimates (or selected predictors) to the "true" parameter or the absolutely best predictor. By absolutely best, we mean the best predictor among those that are linear in $X_1, \ldots, X_m$, rather than the best within a confined subset $B^n$. The goal is to do nearly as well as the absolutely best predictor. In regression when $m \gg n$ there is no hope, in general, to do as well as the absolutely best linear predictor. A natural approach is to confine ourselves to various subsets $B^n$ of the set of all predictors linear in $X_1, \ldots, X_m$, for example, the sets $B^n = A(k)$, where $A(k)$ denotes the set of all the linear predictors which are functions only of $k = k(n)$, $k < m$, explanatory variables. Then we should try to find a predictor which is nearly as good as the corresponding $\beta^*_{F_n}$. Of course, the larger $B^n$, the more challenging is this task. Yet, for too large sets $B^n$, that task is impossible due to reasons explained later using entropy.

It turns out that a sufficient condition for the existence of a persistent sequence of predictors with respect to $B^n$ is that the corresponding sequence $\beta^*_F$ has a sparsity rate $k(n) = o(n/\log(n))$. Note, the last condition on the sparsity rate is trivially satisfied for the sets $B^n = A(k)$, where $k = k(n) = o(n/\log(n))$; hence, our further development is always meaningful for such sets $B^n$.

Our phrasing is slightly different than that of Friedman et al. [11], who write "Use a procedure that does well in sparse problems, since no procedure does well in dense problems." The slight difference in our point of view is that we consider a procedure as doing well, when it does well relative to collection $B^n$ of predictors from which it is feasible to discover nearly the best predictor, with the given sample size. We do not care (since we cannot do much about it) if the absolutely best predictor is indeed in $B^n$ or not. We certainly do not assume that the problem is sparse, that is, that the absolutely best predictor is sparse.

To summarize, we set reasonably high, yet realistic, standards for our procedures, rather than the highest but often impossible to achieve standards.

In Section 2 the procedures achieving persistence will be of the type of best subset selection. More precisely, these procedures search for the empirically best predictor among those in the set $B^n = A(k)$, $k = k(n) = o(n/\log(n))$. Their algorithmic complexity makes such procedures impractical. In Section 3 persistent procedures with lower algorithmic complexity will be introduced for problems with the following intermediate level of generalization. We will consider cases where the function $g_\beta(X_1, \ldots, X_m)$ may be presented as $\rho(\sum \beta_j X_j)$. We will then show that, for those intermediate level of generalization setups, often Lasso-type procedures are useful. By Lasso-type procedures, we mean minimization of $L_{\hat{F}}(\beta)$ subject to a constraint on the $l_1$ norm of $\beta$. Here $\hat{F}$ is the empirical distribution based on the data $Z^1, \ldots, Z^n$ and

$$L_{\hat{F}}(\beta) = \frac{1}{n} \sum_i l(\beta, Z^i).$$

Finally, in Section 4 a simulation study, in high dimensions, is presented for a classification method tied to boosting. The simulated classification method involves optimization under $l_1$ constraint.

*The case where $m < n$.* Our formulation and problems are meaningful also in the case $m < n$. Consider regression again. Let $B^n$ be, as is customary, the set of all linear functions of $X_1, \ldots, X_m$. Then $\beta^*_{F_n}$ is the absolutely best predictor. A related problem in a triangular array formulation was studied by Huber [15], Yohai and Marona [26] and Portnoy [23]; see further references there. In their setup it is desired to estimate the coefficients in a regression problem, where the number of explanatory variables is increased with the number of observations. Under their model, where it is assumed that $Y = \sum \beta_j X_j + \varepsilon$, $E\varepsilon = 0$, the error is not (necessarily) normal and may have heavy tails; also, the explanatory variables are nonrandom. They study consistency in terms of $l_2$ distance between the estimate and the true parameter. The results by Huber and by Yohai and Marona suggest that a sufficient condition for consistency is that the rate that $m$ increases with $n$ is $m = o(\sqrt{n})$. Note that when assuming finite variance for $\varepsilon$, and that the minimal eigenvalue of the design matrix is of order $O(n)$ (as in the case where the columns are orthogonal and the entries are of order 1), a rate $m = o(n)$ is possible using the least squares estimator. However, their interest was mainly in situations involving heavy tails where the variance is not finite.

Portnoy [23] showed that, under natural assumptions, we may let $m$ grow much faster and allow a rate of $m = o(n/\log(n))$. Notice the huge gap compared to the former mentioned rate of $o(\sqrt{n})$. We will also show that the rate suggested by Portnoy should imply persistence in many cases. Yet, we are also left with a similar huge gap; see Remark 4 in Section 2.

**2. Sparsity and persistence.** In this section we will give conditions on triangular arrays under which there exists a procedure satisfying (3).

The following condition will be assumed on the prediction loss $l$ given a triangular array.

CONDITION 1. For every $\varepsilon$, there exists $M(\varepsilon)$, such that for large enough $n$, if $L_{F_n}(\beta) > L_{F_n}(\beta^*_{F_n}) + 2\varepsilon$, then the truncated random variable $T_\varepsilon \equiv \min(l(\beta, Z), M(\varepsilon))$ satisfies

$$E_{F_n} T_\varepsilon > L_{F_n}(\beta^*_{F_n}) + \varepsilon.$$

Note that Condition 1 is obviously satisfied for a bounded prediction loss $l$. In fact, under Condition 1 we may later assume w.l.o.g. that $l(\beta, Z)$ is bounded uniformly under all the distributions in $\mathcal{F}_n$. This will enable us to apply large deviations principles on the fluctuations of $L_{\hat{F}}(\beta)$ from its mean $L_F(\beta)$.

The following easily proved theorem, Theorem 1, is stated for a general triangular array setup. It is a key theorem to understand why, for very general triangular array setups, a predictor should be searched among the set $A(k_n)$ of predictors with corresponding parameters having at most $k_n = o(n/\log(n))$ active entries. In Theorem 6 of [13], it is shown in a regression setup that this rate cannot be improved, that is, an example is given where the sparsity rate is $k_n = O(n/\log(n))$, in which there exists no persistent procedure (of any kind!).

The idea of that proof applies for more general situations, as treated in the current paper. Thus, it seems that, for quite general triangular arrays, when $m = n^\alpha$, $\alpha > 1$, the rate $k_n = o(n/\log(n))$ is also an upper bound for achieving persistence.

$\varepsilon$-*entropy.* We will use the concept of $\varepsilon$-*entropy* of a set of predictors indexed by $\beta$, $\beta \in B$, given a collection of distributions $\mathcal{F}$ and a prediction loss $l$. The definition for it is $\varepsilon$-entropy $\equiv \log(N)$, where $N$ is the minimal number of points, denoted $\beta^1, \ldots, \beta^N$, satisfying that for each $\beta \in B$ there exists a point $\beta^j$ such that, for every $F \in \mathcal{F}$, $|L_F(\beta^j) - L_F(\beta)| < \varepsilon$. A set of such $N$ points will be called an $\varepsilon$-grid.

Note, given any $F$, $F \in \mathcal{F}$, in order to select a predictor whose performance is within $\varepsilon$ of the optimal predictor $\beta^*_F$, it is enough to select the best among an $\varepsilon$-grid of points.

REMARK 2. In order to prove the existence of a persistent procedure with respect to a sequence $B^n$, it is enough to show the existence of a sequence of procedures satisfying (3) for every fixed $\varepsilon$. Then, a diagonalization argument implies the existence of a persistent procedure. Hence, in the following and throughout we will concentrate on showing, for any $\varepsilon > 0$, the existence of a procedure (depending on $\varepsilon$) satisfying (3).

THEOREM 1. *Given a triangular array satisfying Condition 1, assume the following:*

(i) *For every sequence $F_n$, $n = 1, 2, \ldots$, the parameter $\beta^*_{F_n}$ belongs to a $k_n = o(n/\log(n))$-dimensional cube centered at the origin, with Euclidean volume $R_n$, where $\log(R_n) = o(n)$. [Note, in particular, the implied sparsity rate is $o(n/\log(n))$.]*

(ii) *The functions $L_{F_n}(\beta)$ satisfy the following Lipschitz condition: for any $\varepsilon > 0$, there exist $\delta > 0$ and $\gamma > 0$, such that if $\|\beta - \beta'\|_2 < \delta n^{-\gamma}$, $\beta, \beta' \in A(k_n)$, then $|L_{F_n}(\beta) - L_{F_n}(\beta')| < \varepsilon$, uniformly in $F_n \in \mathcal{F}^n$.*

*Then, for every $\varepsilon > 0$ there exists a sequence of procedures satisfying (3), and whence, there exists a persistent procedure.*

CONVENTION. Throughout, we require conditions to hold at $\beta^*_{F_n}$ or at $\hat{\beta} = \arg\min_{\beta \in B^n} L_{\hat{F}}(\beta)$. When these points are not unique, such a condition should be understood as being satisfied if it holds for one of the relevant points.

PROOF OF THEOREM 1. The proof is based on a simple entropy calculation. There are less than $m^{k_n}$ subsets of coordinates of size $k_n$. For each such subset, consider all the predictors determined by active parameters in this subset. For any $F_n \in \mathcal{F}^n$, the function $L_{F_n}(\beta)$, confined to this subset, is viewed as a $k_n$-dimensional function. Divide the corresponding $k_n$-dimensional cube into disjoint small cubes with vertices of length $\delta/\sqrt{k_n} n^\gamma$. Thus, each point in the cube is within Euclidean distance $\delta n^{-\gamma}$ from the center of one of the small cubes, in particular, its true for the point $\beta^*_{F_n}$. These centers determine an $\varepsilon$-grid with respect to the confined versions of $L_{F_n}(\beta), F_n \in \mathcal{F}_n$, given a specific subset and a corresponding $k_n$-dimensional cube; this follows from the Lipschitz condition (ii). The cardinality of the defined $\varepsilon$-grid is $R_n/[\delta/\sqrt{k_n} n^\gamma]^{k_n} = \exp(\log(R_n) + [\log(\frac{1}{\delta}) + \log(\sqrt{k_n}) + \gamma \log(n)]k_n) \equiv A_n$. There are less than $B_n \equiv m^{k_n} = \exp(\alpha \log(n) k_n)$ such subsets, so altogether, the number of points needed to construct an $\varepsilon$-grid, with respect to the set of all predictors containing only points $\beta$ with at most $k_n$ nonzero coordinates and belonging to a cube as in (i), is less than $N = A_n \times B_n$. Now, $\log(A_n \times B_n)$ is of order $o(n)$ if $k_n = o(n/\log(n))$.

It is now standard to show that selecting $\arg\min L_{\hat{F}}(\beta)$, where the minimization is over an $\frac{\varepsilon}{2}$-grid, will yield a procedure that satisfies (3). The reason is as follows: by Condition 1, we may, w.l.o.g., assume that $l$ is bounded and thus, we may conclude exponential rates of convergence to zero of probabilities of large deviations (see, e.g., Hoeffding's inequality [25], page 185). Let $C_n$ be the $\frac{\varepsilon}{2}$-grid of points. Since $\log(N)$ is of order $o(n)$, where $N$ is the cardinality of $C_n$, we obtain, by applying large deviation exponential rates coupled with Bonferroni, $P_{F_n}(\sup_{\beta \in C_n} L_{F_n}(\beta) - L_{\hat{F}_n}(\beta) > \varepsilon) \to 0$. The result now follows. $\square$

The conditions of Theorem 1 imply persistence of procedures that are confined to search for the empirically best predictor among those in $A(k_n)$. The search is also restricted to points which are located in a predetermined cube centered at the origin, where the log of its volume is of order $o(n)$. The best point is "known" to be in such a cube. The last restriction is very weak since the volume of the cube may grow fast. Still, the last restriction could be "mathematically annoying." Under condition (a) of the following Corollary 1, this restriction may be avoided.

Another issue is that the procedure achieving persistence in the proof of Theorem 1 searches in a predetermined grid of points. This is again an artificial restriction. Condition (b) in Corollary 1 requires an analog of the Lipschitz condition in Theorem 1 to hold under the empirical function, $L_{\hat{F}}(\beta)$. Then, it may be concluded that the empirical risk minimization procedure, minimizing over the entire set $A(k_n)$, is persistent, that is, there is no need to minimize in a predetermined set of grid points.

COROLLARY 1. *Consider a triangular array satisfying Condition* 1. *Assume condition* (ii) *of Theorem* 1. *Assume further a sparsity rate* $k_n = o(n/\log(n))$. *Finally, assume the following:*

(a) *With probability approaching* 1 *uniformly for sequences* $F_n$, $\hat{\beta} = \arg\min_{\beta \in A(k_n)} L_{\hat{F}_n}(\beta)$ *belongs to a* $k_n = o(n/\log(n))$-*dimensional cube centered at* $\beta^*_{F_n}$, *with Euclidean volume* $R_n$, *where* $\log(R_n) = o(n)$.

(b) *With probability approaching* 1 *uniformly in sequences* $F_n$, *the random function* $L_{\hat{F}_n}(\beta)$ *satisfies the following Lipschitz condition: for any* $\varepsilon > 0$, *there exist* $\delta > 0$ *and* $\gamma > 0$, *such that if* $\|\beta - \beta'\|_2 < \delta n^{-\gamma}$, *then* $|L_{\hat{F}_n}(\beta) - L_{\hat{F}_n}(\beta')| < \varepsilon$, $\beta, \beta' \in A(k_n)$.

*Then the procedure* $\hat{\beta} = \arg\min_{\beta \in A(k_n)} L_{\hat{F}_n}(\beta)$ *is persistent.*

PROOF. Condition (b) implies that minimizing with respect to $A(k_n)$ is asymptotically equivalent to minimizing with respect to a predetermined (dense enough) grid contained in $A(k_n)$. Similarly, condition (a) implies that minimizing with respect to $A(k_n)$ is asymptotically equivalent to minimizing with respect to its intersection with a predetermined cube centered at $\beta^*_F$. The conclusion now follows by applying Theorem 1. $\square$

Note, often condition (ii) of Theorem 1 follows from condition (b) of Corollary 1.

10 E. GREENSHTEINUnder condition (b) of Corollary 1 and condition (ii) of Theorem 1, for a bounded $l$, the following may be proved similarly to the proof of Corollary 1. Denote by $\mathcal{A}(k_n, R_n)$ the union of all $k_n$-dimensional cubes with volume $R_n$ each. Suppose $\log(R_n) = o(n)$ and $k_n = o(n/\log(n))$. Let $\varepsilon > 0$, and $F_n \in \mathcal{F}_n$ be a sequence of distributions. Then

$$(4) \qquad P_{F_n}\left(\sup_{\beta \in \mathcal{A}(k_n, R_n)} |L_{\hat{F}_n}(\beta) - L_{F_n}(\beta)| > \varepsilon \right) \to 0.$$

*VC dimension.* There is another approach to obtain the type of result in Corollary 1, that is, avoiding the annoying assumption that the optimal predictor is located in a huge cube, and avoiding artificial procedures that search for a predictor in a predetermined grid. It is related to the sophisticated and deep concept of VC dimension.

A way of showing that selecting the predictor that empirically minimizes the risk is equivalent to a search on a grid of $N$ points is through the concept of VC dimension of a class of functions. Using this concept, one may also bound $N$. These bounds depend only on properties of the class of functions $l(\beta, z)$, as functions of $z$ and not of the collection of distributions $\mathcal{F}$ that is involved.

Consider the collection of functions $l(\beta, z) \equiv l_\beta(z)$, $\beta \in B^n$, as functions of $z$. Let us confine ourselves to subsets of functions $l_\beta(z)$ parametrized by $\beta$, whose parameter $\beta$ may have nonzero entries only for certain $k_n$ indices. Suppose the VC dimension of each such confined subset of functions is of order $O(k_n)$. Ideas as in Theorem 1 and Corollary 1 imply that the procedure $\hat{\beta} = \arg\min_{\beta \in A(k_n)} L_{\hat{F}_n}(\beta)$ is persistent when $k_n = o(n/\log(n))$ and a $k_n$ sparsity rate is assumed.

In the following Example 1, we rederive and generalize a result of Greenshtein and Ritov [13]. This is by a simple application of Theorem 1 and Corollary 1. Unlike here, Greenshtein and Ritov used properties of the minimal eigenvalue of a Wishart matrix to establish their result.

EXAMPLE 1. Let $Z_n^i = (Y^i, X_1^i, \ldots, X_{m(n)}^i)$, $m = n^\alpha, \alpha > 1$, where $Z_n^i$ are i.i.d. multivariate normal of dimension $m(n) + 1$ with bounded second moments for $X_j^i$ and $Y^i$ under $\{\mathcal{F}_n\}$. Consider a regression setup, that is, a squared prediction loss $l$ and the set of linear predictors. Under these conditions, we will show that, for $B^n = A(k_n)$, where $k_n = o(n/\log(n))$, the procedure $\hat{\beta} = \arg\min_{\beta \in A(k_n)} L_{\hat{F}_n}(\beta)$ is persistent.

Now, by appropriate reparametrization and invariance considerations, we may assume w.l.o.g., that $X_j$, $j = 1, \ldots, m$, are uncorrelated standard normals; also, w.l.o.g., $Y$ is uncorrelated with the explanatory variables, that is, $\beta_{F_n}^* = 0$. Let $\text{var}(Y) = \sigma^2$. Then $L_F(\beta) = \|\beta\|_2^2 + \sigma^2$, and hence

$$(5) \qquad |L_{F_n}(\beta) - L_{F_n}(\beta')| = |\|\beta\|_2^2 - \|\beta'\|_2^2| \leq \|\beta - \beta'\|_2^2,$$

and condition (ii) of Theorem 1 is satisfied.

In the following we will check conditions (a) and (b) of Corollary 1, in order to finally apply that corollary.

First, the Lipschitz condition, condition (b) of Corollary 1, is satisfied by $L_{\hat{F}_n}(\beta)$ with probability approaching 1. Observe that, for large enough $\gamma'$, $P(\max(X_1,\ldots,X_m) > n^{\gamma'} \equiv M)$ approaches 0, this by combining Chebyshev and Bonferroni. Thus, with high probability, for $\beta, \beta' \in A(k_n)$, $|\sum \beta_j X_j^i - \sum \beta_j' X_j^i| < M \sum |\beta_j - \beta_j'| < M\sqrt{2k_n}\|\beta - \beta'\|_2$. The last inequality is by Cauchy–Schwarz. Condition (b) of Corollary 1 follows, for a squared loss $l$, from the last inequality, when applied similarly to $Z^i$, $i = 1, \ldots, n$.

Condition 1 follows from the multivariate normality. In fact, for the set of random variables $l(\beta, Z)$ with $\|\beta\|_2 < R$, for some $R < \infty$, we have uniform integrability and thus, w.l.o.g., the set consists of bounded random variables.

We now turn to condition (a). We will show that, with probability approaching 1, $\hat{\beta} = \arg\min_{\beta \in A(k_n)} L_{\hat{F}_n}(\beta)$ belongs to a ball with radius (say) $2\sigma^2$, centered at $\beta_{F_n}^* = 0$. Let $G$ be the union of all $k_n$-dimensional balls of radius $2\sigma^2$. Then by the above and by (4), given any $\varepsilon_0 > 0$,

$$(6) \qquad P_{F_n}\left(\sup_{\beta \in G} |L_{F_n}(\beta) - L_{\hat{F}_n}(\beta)| > \varepsilon_0\right) \to 0.$$

Note, that since w.l.o.g. $\beta_{F_n}^* = 0$, we have (*) $L_{F_n}(\beta_{F_n}^*) = L_{F_n}(0) = \sigma^2$. For $\beta$ on the boundary of $G$, $\|\beta\|^2 = 2\sigma^2$, hence, for such $\beta$ we have (**) $L_{F_n}(\beta) = 2\sigma^2 + \sigma^2$. Condition (a) now follows by the convexity of $L_{\hat{F}_n}(\beta)$, from (*), (**) and (6).

Finally, applying Corollary 1, we obtain that the procedure $\hat{\beta}$ which selects the empirically best predictor from the set $A(k_n)$ is persistent.

REMARK 3. (i) In the last example we used only multivariate normality to conclude Condition 1. Hence, the result holds in much more general situations. A proof along the lines of Example 1 is possible for other prediction losses, for example, $l(\beta, Z) = |Y - \sum \beta_j X_j|$.

REMARK 4. Consider a regression case, as in Example 1. Suppose we replace the multivariate normal assumption by the assumption that the entries of $Z^i$ are bounded under the possible distributions in the triangular array. We cannot prove the $o(n/\log(n))$ rate for $k_n$, as in Example 1. The reason is that Condition 1 is not implied. Note, existence of $M(\varepsilon)$ for every fixed $n$ is trivially implied by boundedness, but not existence of $M(\varepsilon)$ that holds uniformly for every $n$. In [13] a sparsity rate of $k_n = o(\sqrt{n/\log(n)})$ is shown to imply persistence, under an additional assumption, that the minimal eigenvalue of the covariance matrix of $(X_1, \ldots, X_m)$ does not approach

0. Whether we may obtain persistence under higher rates, assuming only boundedness, is suggested there as a problem. We still do not know the answer to this problem.

**3. Optimization under $l_1$ constraint and the Lasso.** In the main result of this section, Theorem 2, we will show for special classes of parametrized predictors that we may achieve persistence and approximate the best subset of a certain size through optimization under $l_1$ constraint. The special classes are of the form

$$(7) \qquad g_\beta(X_1, \ldots, X_m) = \rho\Big(\sum \beta_j X_j\Big).$$

As a further example, consider the class of predictors

$$(8) \qquad g_\beta(X_1, \ldots, X_m) = \frac{\exp[\sum \beta_j X_j]}{1 + \exp[\sum \beta_j X_j]}.$$

The optimization under the constraint that the number of nonzero entries of $\beta$ is $k_n$ has high complexity in general. It is desired to replace it by a constraint that determines a convex feasible set. When the target function $L_{\hat{F}}(\beta)$ is also convex, then the problem has an algorithmically efficient solution; see [20].

An example where both the target function and the feasible set are convex is the Lasso procedure, that is,

$$(9) \qquad \min_\beta L_{\hat{F}}(\beta) = \min_\beta \frac{1}{n} \sum_i \Big(Y^i - \sum \beta_j X_j^i\Big)^2,$$

subject to the constraint $\|\beta\|_1 < b$ for a proper $b$. See [24]; also see basis pursuit in [5]. Recently Efron et al. [8] developed an efficient algorithm, called least angle regression, to solve the above optimization problem. We will elaborate on another example involving convex optimization in the next section.

We study the replacement of the constraint on the number of nonzero entries of $\beta$ by a convex constraint on its $l_1$ norm. In recent papers by Donoho [6] and [7], a general setup is described, in which optimization under $l_1$ constraint gives the actual optimal solution under the constraint on the number of nonzero entries. Our ultimate goal is not to find a predictor with a sparse representation; for us, searching for a sparse solution is only a means of regularization and of controlling the entropy. Thus, we need weaker results compared to those of Donoho; for our purpose, it is enough to show some kind of (weaker) equivalence between the solutions obtained under the two types of constraints. From the following Lemma 1, it follows that predictors with parameters that are obtained through optimization under a constraint on their $l_1$ norm might (appear to) have more than $k_n$ "active entries," but in fact it will be shown that, keeping the $l_1$ constraint in the right magnitude (depending on $k_n$), they are equivalent to predictors with parameters that have only $k_n$ active entries. From the last fact, our main theorem of this section, Theorem 2, will follow. It is a generalization of a result obtained by Greenshtein and Ritov [13] for regression.

The following Lemma 1 is given without proof. It is a rephrasing of a result by Maurey (see [22]; a version of it may be found in [13], Lemma 4, and in [16], Proposition 2.2). There, the analogous result is stated for a single distribution $G$, but the same proof works for a pair $G_1$ and $G_2$, as in what follows.

LEMMA 1. *Let $G_1$ and $G_2$ be two distributions under which $X_j$, $j = 1, \ldots, m$, are bounded by $M$. Let $\beta$ be an $m$-dimensional vector such that $\|\beta\|_1 = b$. Let $\delta > 0$. Then for every $\kappa > 0$, there exists a corresponding vector $\beta'$, where $\|\beta'\|_1 = b$, having at most $\kappa$ nonzero coefficients, such that*

$$P_{G_i}\Big(\Big|\sum \beta_j X_j - \sum \beta'_j X_j\Big| > \delta\Big) < M^2 b^2/\delta^2 \kappa, \qquad i = 1, 2.$$

We will confine ourselves to triangular arrays where, for each $n$, the pair consisting of prediction loss $l$ and the collection of predictors $\{g_\beta\}$ satisfies the following:

CONDITION 2. *For a fixed $y$, the function*

$$h\Big(y, \sum \beta_j X_j\Big) \equiv l(Y, g_\beta(X_1, \ldots, X_m))$$

*is bounded and uniformly continuous in $\sum \beta_j X_j$, uniformly in $y$.*

The boundedness condition on $l$ may be circumvented in various examples. It may be weakened assuming a condition like Condition 1, or uniform integrability of $l(\beta, Z)$, $\beta \in B^n$. In Theorem 2 we will also require boundedness of $X_j^i$, this is in order to apply Lemma 1. If this assumption is avoided, the required sparsity rate in Theorem 2 would be $o(n/\log(n)d_n)$, where $d_n = \sup_{\mathcal{F}_n} E_{F_n}[\max(X_1, \ldots, X_m)]^2$. Again, the boundedness assumption on $X_j^i$ may be avoided in special cases, like regression with multivariate normal $Z^i$, as treated in Section 4 of [13]. We will leave the boundedness assumption for a clearer exposition.

Our main theorem for this section is the following.

THEOREM 2. *Consider a triangular array satisfying Condition 2 and having bounded $X_j^i$. Suppose the sparsity rate is $k_n = o(n/\log(n))$. Suppose further that $\|\beta_{F_n}^*\|_2$ is bounded by $R$ (w.l.o.g. $R = 1$) for every $F_n$, $F_n \in \mathcal{F}_n$,* $n = 1, 2, \ldots$. *Then the following procedure is persistent. Select the predictor* $\tilde{\beta}$, *where*

$$\tilde{\beta} = \arg\min_{\beta} L_{\hat{F}_n}(\beta), \tag{10}$$

*subject to the constraints* $\|\beta\|_1 \leq \sqrt{k_n}$.

LEMMA 2. *Assume $X_j^i$ are bounded by $M$. Then Condition 2 implies both Lipschitz conditions, that is, condition* (ii) *in Theorem 1 and condition* (b) *in Corollary 1.*

PROOF. Observe that

$$\left|\sum \beta_j X_j - \sum \beta_j^* X_j\right| < M \sum |\beta_j - \beta_j^*| \tag{11}$$
$$< M\sqrt{k_n}\|\beta - \beta^*\|_2 < Mn^{0.5}\|\beta - \beta^*\|_2.$$

Here we have applied Cauchy–Schwarz and the fact that $k_n < n$.

The proof follows from the uniform continuity and boundedness of $l$. □

PROOF OF THEOREM 2. Let $\tilde{\beta}$ be the solution of (10) for a data set coming from $F_n$. Then by Lemma 1, given $\varepsilon_1 > 0$ and $\delta_1 > 0$, for any sequence $\kappa_n$ such that $k_n = o(\kappa_n)$, there exists a parameter $\beta'$ having at most $\kappa_n$ nonzero entries, such that both for $G_1 = F_n$ and for $G_2 = \hat{F}_n$ we have

$$P_{G_i}\left(\left|\sum \beta_j' X_j - \sum \tilde{\beta}_j X_j\right| > \delta_1\right) < \varepsilon_1, \qquad i = 1, 2. \tag{12}$$

Moreover,

$$\|\beta'\|_1 = \|\tilde{\beta}\|_1 \leq \sqrt{k_n}. \tag{13}$$

We choose a sequence $\kappa_n$ which is $o(n/\log(n))$, so that (12) is satisfied.

By Condition 2, (12) implies both

$$|L_{F_n}(\tilde{\beta}) - L_{F_n}(\beta')| < \varepsilon = \varepsilon(\varepsilon_1, \delta_1) \tag{14}$$

and

$$|L_{\hat{F}_n}(\tilde{\beta}) - L_{\hat{F}_n}(\beta')| < \varepsilon = \varepsilon(\varepsilon_1, \delta_1). \tag{15}$$

Note that we may obtain (14) and (15) for $\varepsilon > 0$ arbitrarily small, by selecting large enough $\kappa_n = o(n/\log(n))$. By (13) and by construction, $\beta'$ belongs to a $\kappa_n = o(n/\log(n))$-dimensional cube centered at $\beta^*_{F_n}$, where the logarithm of the cube's volume is $o(n)$. Also, by Lemma 2, both condition (b) of Corollary 1 and condition (ii) of Theorem 1 are satisfied. Hence, by (4) we obtain

$$P_{F_n}(|L_{\hat{F}_n}(\beta') - L_{F_n}(\beta')| > \varepsilon) \to 0. \tag{16}$$

Note, by assumption $\|\beta^*_{F_n}\|_2 \leq 1$, whence, by Cauchy–Schwarz, $\|\beta^*_{F_n}\|_1 \leq \sqrt{k_n}$. Thus, by the definition of $\tilde{\beta}$ we have

(17) $$L_{\hat{F}_n}(\tilde{\beta}) \leq L_{\hat{F}_n}(\beta^*_{F_n}).$$

Finally, by the law of large numbers we have

(18) $$P_{F_n}(|L_{F_n}(\beta^*_{F_n}) - L_{\hat{F}_n}(\beta^*_{F_n})| > \varepsilon) \to 0.$$

From (14), (15), (16), (17) and (18), we obtain persistence of $\tilde{\beta}$, and the proof of the theorem follows. □

As remarked before, in practice the proper value for the $l_1$ constraint is unknown. One should try various values and test the resulting predictors on a test set. Our theory suggests that the resulting optimal $l_1$ constraint will be of order $\sqrt{n/\log(n)}$.

REMARK 5. From the proof of Theorem 2, we obtain, even when not assuming sparsity, an appealing feature of rules based on optimization under $l_1$ constraint. The feature is *self consistency* of such procedures. The self consistency is in the following sense.

Suppose $\tilde{\beta}$ is obtained by (10) for $k_n = o(n/\log(n))$. Suppose Conditions 1 and 2 hold. Then for every $\varepsilon > 0$ and every sequence $F_n$,

$$P_{F_n}(|L_{\hat{F}_n}(\tilde{\beta}) - L_{F_n}(\tilde{\beta})| > \varepsilon) \to 0.$$

COROLLARY 2. *From the above it follows that, under Condition* 2, *the procedure defined by* (10) *is persistent with respect to* $B^n$, *the sequence of $l_1$ balls with an $l_1$ radius of order $k_n = o(\sqrt{\frac{n}{\log(n)}})$.*

*(There is no need to assume sparsity.)*

*Discussion. Regularization by general $l_q$ constraints.* The $l_1$ constraint is motivated through a constraint on the number of nonzero parameters, which may also be represented as an $l_0$ constraint. The advantage of the $l_1$ constraint relative to other $l_q$ constraints, $q < 1$, is the convexity of the feasible set. Yet, from Theorem 2, we conclude that we will not gain much by optimizing via an $l_0$ or $l_q, q < 1$, constraint. This is since persistence under a $o(n/\log(n))$ sparsity rate is already achieved using $l_1$ constraint, while the proofs in this paper and the forementioned Theorem 6 of Greenshtein and Ritov [13] indicate that, in general, persistence cannot be achieved for higher rates.

*Lack of persistence of ridge regression.* Regularization via $l_q$ constraint with $q > 1$ will usually lead to nonpersistent procedures, which are also not self-consistent. Consider, for example, the case $q = 2$ in a regression context with a squared loss, called ridge regression. Suppose $Z^i$, $i = 1, \ldots, n$, are multivariate normal, and suppose $\beta^*_{F_n} = 0$, that is, $Y^i$ are not correlated with the corresponding $m$ explanatory variables. Assume also that $X_j$ are uncorrelated standard normals. Denote $\sigma^2 = \text{var}(Y)$. Minimizing the empirical risk subject to a constraint $\sum \beta_j^2 < \delta^2$ will yield (typically) a solution $\hat{\hat{\beta}}$ which is on the boundary of the feasible set when $\delta^2 < \sigma^2$, that is, the estimate will have an $l_2$ norm $\delta$. This situation remains when $m$ and $n$ approach infinity in a way that $m \gg n$, no matter how small is $\delta > 0$. Thus, $L_{F_n}(\hat{\hat{\beta}}) = \delta^2 + \sigma^2 + o_p(1)$, while $L_{F_n}(\beta^*_{F_n}) = L_{F_n}(0) = \sigma^2$.

When the regularization is via an $l_1$ constraint, as suggested in this paper, again the minimizer of the empirical risk, denoted $\hat{\beta}$, will be on the boundary of the $l_1$ ball, which is the feasible set. Yet now, when the $l_1$ constraint is chosen properly to be $o(\sqrt{n/\log(n)})$, the $l_2$ norm of that solution will be of order $o_p(1)$, hence, $L_{F_n}(\hat{\beta}) = \sigma^2 + o_p(1)$. This property of the $l_1$ constraint is a consequence of our Theorem 2.

Further discussion of the $l_1$ constraint regularization method and its comparison with $l_2$ regularization may be found in [11] and [4].

In general, regularization may be achieved by introducing penalty functions. For example, using Lagrange multipliers, one may see that the solution of the optimization problem, under $l_q$ constraint, is the same as the solution of the related optimization problem when introducing the penalty function $\lambda \sum |\beta_j|^q$, called $l_q$ penalization, for an appropriate $\lambda$. A study of regularization using general penalizations was conducted by Fan and Li [9] and by Fan and Peng [10]. In their setup analogous to our prediction loss $l(\cdot)$ is the log-likelihood, but the essence is the same (see some elaboration on it in [12]). They treat a general class of penalty functions, including the $l_q$ penalties. In particular, for $l_q$ penalization with $q < 1$ and a proper choice of $\lambda$, they show that a certain oracle optimality is achieved by penalized maximum likelihood procedures, while for $q = 1$, such optimality does not seem to be implied (the recommended penalty functions in those papers are not an $l_q$ type, but a class of penalty functions called SCAD which possesses further nice properties). In a sparse setup, an oracle optimality of procedures means the following. The rate of convergence to the estimated parameter is the same as the rate that may be achieved when knowing which are the zero entries of the parameter. These results are obtained also under a triangular array setup in [10], but when $m(n) \ll n$. In particular, for $m = o(n^\alpha)$, $\alpha = \frac{1}{5}, \frac{1}{4}, \frac{1}{3}$, under various assumptions and regularity conditions. These oracle optimality properties are much more delicate and strong than the persistence suggested by us. Such strong optimality criteria may be achieved by procedures, due to the slow rate at which the dimension $m = m(n)$ increases with $n$, in comparison to the rate in our setting.

**4. Numerical study.** In this section we examine through simulation the following high-dimensional classification problem. Consider $Z = (Y, X_1, \ldots, X_m)$, where the value of $Y$ is either $-1$ or $+1$. The prediction loss is

(19) $$l(\beta, Z) = h\Big(Y, \sum \beta_j X_j\Big) = \exp\Big(-Y \sum \beta_j X_j\Big).$$

The convex loss (19) is used to motivate the boosting classification procedure; see, for example, [14], page 305, or [3]. It may also be motivated as follows. Suppose we classify according to $g_\beta(X_1, \ldots, X_m) = \text{sign}(\sum \beta_j X_j)$. Now the value of $\sum \beta_j X_j$ is interpreted both through its sign and the magnitude of its absolute value. The sign determines the classification decision and the magnitude is interpreted as the "confidence in that decision." That is why wrong classifications with large magnitude are severely penalized and vice versa.

Our optimization under $l_1$ constraint is similar to the approach of Lugosi and Vayatis [18]. As observed by them, there could be many other interesting and natural convex prediction losses other than the above; for example, see their Example 3. Yet, (19) has attracted a lot of attention recently and we elaborate on it.

In the following we present a simulation study where the dimension $m$ is of the order of thousands, while the sample size $n$ is of the order of hundreds.

*The simulation.* We simulate $n$ i.i.d. vectors. Each is $M$-dimensional and consists of $M$ i.i.d. $N(0,1)$, random variables. Denote the $j$th component of the $i$th vector by $X_j^i$.

For each vector $i$, $i = 1, \ldots, n$, let $W^i$ be a $N(0, 0.25)$ random number independent of $X_j^i$ and define

$$Y^i = \text{sign}\bigg(\frac{X_1^i + \cdots + X_{25}^i}{5} + W^i\bigg) = \text{sign}(V^i + W^i),$$

where $V^i$ is implicitly defined. Thus, the first 25 "explanatory variables" (out of the $M$ available ones) are the relevant predictors for $Y^i$, and the prediction should be through $V^i$.

Now we create, for each $i$, five additional random numbers (or simulated explanatory variables), denoted $X_{M+1}^i, \ldots, X_{M+5}^i$, as follows: $X_j^i = V^i + U_j^i$, $j = M+1, \ldots, M+5$; here $U_j^i \sim N(0, 9)$ are again independent of all the others and of each other.

Notice we have $m = M + 5$ explanatory variables; only the first 25 are relevant for predicting $Y^i$. Yet, if we may choose only a single explanatory variable to base our prediction on, we would rather choose $X_j^i$ from the group of the last five; obviously if we may choose as many as 25 or more, we would choose the first 25.

Table 1
$n = 500$ and $M = 1000$ ($m = 1005$)

| V-training | V-real | B1-1-norm | B2-1-norm | $\beta$-1-norm | $\lambda$ |
|---|---|---|---|---|---|
| 0.132 | 2.365 | 7.409 | 0.444 | 22.422 | 0.01 |
| 0.361 | 0.850 | 3.985 | 0.291 | 9.757 | 0.03 |
| 0.538 | 0.810 | 2.277 | 0.270 | 5.030 | 0.05 |
| 0.673 | 0.817 | 1.430 | 0.240 | 2.767 | 0.07 |
| 0.742 | 0.850 | 0.825 | 0.277 | 1.540 | 0.09 |
| 0.815 | 0.860 | 0.499 | 0.246 | 0.887 | 0.11 |
| 0.859 | 0.880 | 0.243 | 0.242 | 0.523 | 0.13 |
| 0.877 | 0.895 | 0.142 | 0.229 | 0.379 | 0.15 |
| 0.887 | 0.902 | 0.084 | 0.224 | 0.311 | 0.17 |

Our indirect method of searching for the best subset is through optimization under $l_1$ constraint. Practically, the right constraint may be determined by cross-validation or a test set. In our simulation study, the performance of a predictor, obtained through such optimization under $l_1$ constraint, was tested on an independent sample of size 1000. In Tables 1–3 the average prediction loss on the "data set"/"training set" is denoted $V$-training, while the average on the additional independent sample of size 1000 is denoted $V$-real.

Our optimization is conducted using "Lagrange multipliers," that is, instead of optimization under $l_1$ constraint, we optimize, for appropriate $\lambda > 0$,

$$L_{\hat{F}}(\beta) + \lambda \sum |\beta_j|.$$

We try various values of $\lambda$ that correspond to various constraints on the $l_1$ norm of $\beta$. The optimization is through steepest descent, where special care is taken when computing the "partial derivative" of $\lambda \sum |\beta_j|$, for coordinates $j$ where for the current iteration $\beta_j = 0$.

In Tables 1–3 we summarize simulation results for various $m$ and $n$. Only for the case $n = 500$, $M = 1000$ is a detailed table given, with the performance under various constraints. For the other cases, $n = 100$, $M = 1000$ and $n = 500$, $M = 5000$, only the performance under the optimal constraint is given. Each row is based on averages of 20 repetitions for a fixed $\lambda$. In the same table, different rows correspond to different $\lambda$, and the bigger $\lambda$ is, the more severe is the constraint. Indeed, one may see in Table 1 that as $\lambda$ decreases the difference between $V$-training and $V$-real increases, that is, the generalization power (or self consistency property) is reduced. We record the constraint also in terms of the $l_1$ norm of $\beta$ in the column $\beta$-1-norm. The columns B1-1-norm and B2-1-norm record the $l_1$ norm of the first 25 and of the last five coordinates, respectively.

Table 2
$n = 100$ and $M = 1000$ ($m = 1005$)

| V-training | V-real | B1-1-norm | B2-1-norm | $\beta$-1-norm | $\lambda$ |
|---|---|---|---|---|---|
| 0.861 | 0.926 | 0.010 | 0.207 | 0.264 | 0.30 |

Table 3
$n = 500$ and $M = 5000$ $(m = 5005)$

| V-training | V-real | B1-1-norm | B2-1-norm | $\beta$-1-norm | $\lambda$ |
|---|---|---|---|---|---|
| 0.690 | 0.862 | 0.680 | 0.271 | 2.181 | 0.09 |

In practice, the column $V$-real will be replaced by evaluation of the performance of the suggested predictor on a test set or cross-validation (the evaluation would be less accurate when the test set is smaller than the 1000 used in our simulation). Thinking of the $V$-real column as results from a test set, we get the following. When there are only $n = 100$ observations available, a test set would suggest to predict mainly based on the last five explanatory variables using $\lambda \approx 0.3$ and with risk $\approx V$-real $= 0.926$. Note, the $l_1$ mass of the first 25 coefficients is only 0.01, while the $l_1$ mass of the last five is 0.207. When there are $n = 500$ observations, a test set would suggest $\lambda \approx 0.05$ with resulting risk about 0.81. Note, when $n = 500$, the $l_1$ mass of the first 25 coefficients is 2.277, while that of the last five is only 0.27. Indeed, with only 100 observations, the attempt to reveal the 25 "best" explanatory variables is too ambitious and the procedure gives up on it and settles for the inferior group of five. When the sample size is increased to 500, there is a shift toward the first 25 variables.

Comparing the simulated results with $M = 1000$ to those with $M = 5000$, we see that by screening in advance many superfluous explanatory variables, reducing from $m = 5005$ to $m = 1005$, we hardly improve. In the case $m = 5005$ the best value is attained when $\lambda = 0.09$ and equals $V$-real $= 0.862$; in the case $m = 1005$ the best value is attained when $\lambda = 0.05$ and equals $V$-real $= 0.810$. The improvement is by 0.052. One could argue that this improvement might be significant when compared to the risk magnitudes, 0.810 and 0.862. As remarked in the Introduction, when the risk is small (or approaches 0), a more delicate analysis of rates of convergence, rather than only persistence, is desired.

Note, however, that the slight advantage demonstrated when screening out successfully 4000 superfluous explanatory variables (in our simulation changing $m$ from 5005 to 1005) seems to occur in the "twilight zone," that is, the zone where the constraint is not severe enough to produce estimators with generalization power (or that are self-consistent). Compare in Table 1 for the optimal constraint $\lambda = 0.05$, $V$-real $= 0.81$, while $V$-training $= 0.538$. Such a "twilight zone" could be very abrupt in very high dimensions. Moving further from that zone will introduce singularity and the selected predictors will be totally unreliable.

**Acknowledgment.** I am grateful to Anirban DasGupta for comments that led to a better presentation.

DEPARTMENT OF STATISTICS
PURDUE UNIVERSITY
WEST LAFAYETTE, INDIANA 47907
USA
E-MAIL: egreensh@stat.purdue.edu